\documentclass[titlepage,12pt]{article}
\usepackage{amssymb}
\usepackage{amsfonts}
\usepackage{amsmath}
\begin{document}

{\bf \Large Archimedean Type Conditions in Categories} \\ \\

{\bf Elemer E Rosinger} \\
Department of Mathematics \\
and Applied Mathematics \\
University of Pretoria \\
Pretoria \\
0002 South Africa \\
eerosinger@hotmail.com \\ \\

{\bf Abstract} \\

Two concepts of being Archimedean are defined for arbitrary categories. \\ \\

{\bf 1. The case of usual semigroups} \\

For convenience, let us recall two versions of concepts of being Archimedean in the usual
case of algebraic structures. In this regard, a sufficiently general setup is as follows. \\

Let $( E, +, \leq )$ be a partially ordered semigroup, thus we have satisfied \\

(1.1)~~~ $ x, y \in E_+ ~~\Longrightarrow~~ x + y \in E_+ $ \\

where $E_+ = \{ x \in E ~|~ x \geq 0 \}$. \\

A first intuitive version of the Archimedean condition, suggested in case $\leq$ is a {\it
linear} or {\it total} order on $E$, is \\

(1.2)~~~ $ \exists~~ u \in E_+ ~:~ \forall~~ x \in E ~:~
                            \exists~~ n \in \mathbb{N} ~:~ n u \geq x $ \\

Here is another formulation used in the literature when $\leq$ is an arbitrary partial order
on $E$, namely \\

(1.3)~~~ $ \begin{array}{l}
                \forall~~ x \in E_+ ~:~ \\
                 ~~~~ x = 0 ~\Longleftrightarrow~
                          \left ( \begin{array}{l}
                                       \exists~~ y \in E_+ ~: \\
                                       \forall~~ n \in \mathbb{N} ~: \\
                                        ~~~~ n x \leq y
                                   \end{array} \right )
          \end{array} $ \\ \\

where clearly the implication "$\Longrightarrow$" is trivial, and which condition is
thus equivalent with \\

(1.4)~~~ $ \forall~~ x \in E_+ ~:~
               \mathbb{N} x ~~\mbox{is bounded above}~~
                  \Longrightarrow ~~ x = 0 $ \\

{\bf Lemma 1.} \\

If $( E, +, \leq )$ is a {\it linearly} or {\it totally} ordered semigroup, then (1.3)
~~$\Longrightarrow$~~ (1.2). \\

{\bf Proof.} \\

Assume indeed that (1.2) does not hold, then \\

$~~~~~~ \forall~~ x \in E_+ ~:~ \exists~~ y \in E ~:~
                            \forall~~ n \in \mathbb{N} ~:~ n x \ngeqslant y $ \\

and since $\leq$ is a linear or total order on $E$, we have \\

$~~~~~~ \forall~~ x \in E_+ ~:~ \exists~~ y \in E ~:~
                            \forall~~ n \in \mathbb{N} ~:~ n x \leq y $ \\

Obviously, we can assume that $y \in E_+$, thus (1.4) is contradicted. \\ \\

{\bf 2. The case of categories} \\

Let ${\cal C}$ be any category. The issue is to be able to take a so called "unit" morphism,
like for instance $u$ in (1.2), say \\

$~~~ A \stackrel{f}\longrightarrow B $ \\

and be able to "repeat" it, say, to the right of $B$ any finite number of time. Here the
problem is that, in general, we cannot compose a morphism in ${\cal C}$ with itself even just
twice. In particular, we cannot in general have $f \circ f$, let alone $f \circ f\circ f$, and
so on. Therefore, when given two ${\cal C}$ morphisms which can be composed \\

$~~~ A \stackrel{f}\longrightarrow B \stackrel{g}\longrightarrow C $ \\

we have to find a way to be able to say that the morphism $g$ is again a "unit", that is, more
or less the same with the morphism $f$ from a certain relevant point of view. \\

One simple natural way to do that is as follows. We consider the {\it arrow category} ${\cal
C}^2$, associated with ${\cal C}$, [H \& S, p. 27], namely, the category whose class of {\it
objects} is the class of ${\cal C}$ morphisms, while for any two such ${\cal C}^2$ objects \\

$~~~ A \stackrel{f}\longrightarrow B,~~~ A\,' \stackrel{f\,'}\longrightarrow B\,' $ \\

the corresponding ${\cal C}^2$ morphisms are the pairs $( a, b )$, where \\

$~~~ A \stackrel{a}\longrightarrow A\,',~~~ B \stackrel{b}\longrightarrow B\,' $ \\

are ${\cal C}$ morphisms, such that the diagram commutes \\

\begin{math}
\setlength{\unitlength}{0.2cm}
\thicklines
\begin{picture}(60,15)

\put(14,12){$A$}
\put(36,12){$B$}
\put(25,13.7){$f$}
\put(17,12.6){\vector(1,0){18}}
\put(13.5,1){$A\,'$}
\put(12.5,7){$a$}
\put(14.7,10.5){\vector(0,-1){6.5}}
\put(37.8,7){$b$}
\put(37,10.5){\vector(0,-1){6.5}}
\put(36,1){$B\,'$}
\put(17,1.6){\vector(1,0){18}}
\put(25,-1){$f\,'$}
\end{picture}
\end{math} \\

Finally, the ${\cal C}^2$ composition of morphisms is defined by \\

$~~~ ( a\,', b\,' ) \circ ( a, b ) = ( a\,' \circ a, b\,' \circ b) $ \\

in other words, by pasting the two above kind of diagrams together, and deleting the middle
horizontal arrow. \\

Now, given two ${\cal C}$ morphisms \\

$~~~ A \stackrel{f}\longrightarrow B,~~~ C \stackrel{g}\longrightarrow D $ \\

we say that they are {\it unitary equivalent}, if and only if, when considered as objects in
the arrow category ${\cal C}^2$, they are isomorphic, [H \& S, p. 36]. \\

With that definition, we can now attempt to define when a category ${\cal C}$ is {\it
Archimedean}. \\

First, we consider an extension of the usual version of the Archimedean property in (1.2).
Namely, the corresponding condition in categories is as follows \\

(2.1) $~~~ \begin{array}{l}
       \exists~~ U \stackrel{v}\longrightarrow W~~~ {\cal C} ~\mbox{morphism} ~: \\ \\
       \forall~~ A \stackrel{f}\longrightarrow B~~~ {\cal C} ~\mbox{morphism} ~: \\ \\
       \exists~~ U_1 \stackrel{v_1}\longrightarrow W_1, \ldots, U_n
                    \stackrel{v_n}\longrightarrow W_n~~~ {\cal C} ~\mbox{morphisms} ~: \\ \\
       ~~~~1)~~~ v_1, \ldots, v_n ~~\mbox{are in that order composable in}~~ {\cal C} \\ \\
       ~~~~2)~~~ v_1, \ldots , v_n ~~\mbox{are unitary equivalent with}~~ v \\ \\
       ~~~~3)~~~ f ~~\mbox{is a submorphism of}~~ v_n \circ \ldots \circ v_1
            \end{array} $ \\ \\

Here we used the following definition. Given two ${\cal C}$ morphisms \\

$~~~ A \stackrel{f}\longrightarrow B,~~~ C \stackrel{g}\longrightarrow D $ \\

we say that $f$ is a {\it submorphism} of $g$, if and only if there are ${\cal C}$
morphisms \\

$~~~ C \stackrel{a}\longrightarrow A,~~~ B \stackrel{b}\longrightarrow D $ \\

such that \\

(2.2) $~~~ g = b \circ f \circ a $ \\

As for the extension to categories of the usual condition (1.4), we can proceed as follows.
Given any ${\cal C}$ morphism $U \stackrel{v}\longrightarrow W$, let us denote by \\

(2.3) $~~~ \mathbb{N} v $ \\

the class of ${\cal C}$ morphisms of the form \\

(2.4) $~~~ v_n \circ \ldots \circ v_1 $ \\

where $U_1 \stackrel{v_1}\longrightarrow W_1, \ldots, U_n \stackrel{v_n}\longrightarrow W_n$
are ${\cal C}$ morphisms which are composable in the given order, and are each unitary
equivalent with $v$. \\

Then we call ${\cal C}$ {\it Archimedean}, if and only if \\

(2.5) $~~~ \begin{array}{l}
                \forall~~ U \stackrel{v}\longrightarrow W~~~
                                  ~ {\cal C} ~\mbox{morphism} ~: \\ \\
                ~~~~ \mathbb{N} v ~~\mbox{ is bounded in}~~ {\cal C}
                         ~~~\Longrightarrow~~~ v = id_U
            \end{array} $ \\ \\

Here the following definition was used. A given class ${\cal N}$ of ${\cal C}$ morphisms is
called {\it bounded}, if and only if there exists a ${\cal C}$ morphism $A \stackrel{f}
\longrightarrow B$, such that every ${\cal C}$ morphism in ${\cal N}$ is a submorphism of
$f$. \\ \\

{\bf 3. Example} \\

As a simple example we shall illustrate the two general concepts of Archimedean category
presented in section 2, in the particular case of categories given by {\it quasi ordered}
classes, [H \& S, p. 19]. We recall that a category ${\cal C}$ is a quasi ordered class, if
and only if for every two of its objects $A$ and $B$, there is at most one single morphism $A
\stackrel{f}\longrightarrow B$. It follows that, for simplicity, such a category can be
described by a {\it reflexive} and {\it transitive} binary relation $\leq$ on its objects,
which is defined as follows. For every two objects $A$ and $B$ in ${\cal C}$, we have $A \leq
B$, if and only if there exists a morphism $A \stackrel{f}\longrightarrow B$ in ${\cal C}$. In
such a case that morphism is, as assumed, unique. \\

In order to illustrate in the above particular situation the two concepts of being Archimedean
defined in section 2 for arbitrary categories, we first have to clarify in the context of
quasi ordered classes the notion of {\it unitary equivalent}, which plays a role in both
mentioned concepts. And for this purpose, we have to look at the arrow category ${\cal C}^2$ of
the quasi ordered class ${\cal C}$. \\

It is easy to see that the objects in ${\cal C}^2$ are precisely the pairs $( A, B )$, where
$A$ and $B$ are objects in ${\cal C}$, such that $A \leq B$. Further, given two objects $( A,
B )$ and $( A\,', B\,' )$ in ${\cal C}^2$, it is immediate that, in  ${\cal C}^2$, there exist
morphisms between them, if and only if \\

(3.1) $~~~ A \leq A\,',~~~ B \leq B\,' $ \\

And in such a case, the corresponding morphism is the pair of pairs \\

(3.2) $~~~ ( ( A, A\,' ), ( B, B\,' ) ) $ \\

which is the unique such morphism in ${\cal C}^2$. \\

Thus ${\cal C}^2$ is again a quasi ordered class. \\

As for the composition of such morphisms in ${\cal C}^2$, one readily obtains that \\

(3.3) $~~~( ( A, A\,' ), ( B, B\,' ) ) \circ ( ( A\,', A\,'' ), ( B\,', B\,'' ) ) =
                      ( ( A, A\,'' ), ( B, B\,'' ) ) $ \\

Given now two morphisms $( A, B )$ and $( C, D )$ in ${\cal C}$, then by the definition in
section 2, they are unitary equivalent, if and only if, when considered as objects in ${\cal
C}^2$, they are isomorphic. In other words, there must be an isomorphism $( A, B )
\stackrel{\alpha}\longrightarrow ( C, D )$ in ${\cal C}^2$. This means, [H \& S, Proposition
5.17, p. 36], that there exists a unique morphism $( C, D ) \stackrel{\beta}\longrightarrow
( A, B )$ in ${\cal C}^2$, such that $\beta \circ \alpha = id_{( A, B )}$ and $\alpha \circ
\beta = id_{( C, D)}$. \\
Thus in view of (3.2), it follows that \\

$~~~ \alpha = ( ( A, C ), ( B, D ) ),~~~ \beta = ( ( C, A ), ( D, B ) ) $ \\

Furthermore, (3.1) applied to $\alpha$, implies \\

(3.4) $~~~ A \leq C,~~~ B \leq D $ \\

and similarly, applied to $\beta$, gives \\

(3.5) $~~~ C \leq A,~~~ D \leq B $ \\

Now as customary, let us consider the {\it equivalence} relation $\approx$ on the objects of
the quasi ordered class ${\cal C}$, defined by \\

(3.6) $~~~ X \approx Y ~~~\Leftrightarrow~~~ X \leq Y,~~~ Y \leq X $ \\

Then (3.4), (3.5) give for any two morphisms $( A, B )$ and $( C, D )$ in ${\cal C}$ \\

(3.7) $~~~ \begin{array}{l}
             (  A, B ),~( C, D ) ~~\mbox{are unitary equivalent} ~~~\Longleftrightarrow \\ \\
                 ~~~~~~~~~~~\Longleftrightarrow~~~ A \approx C \leq B \approx D
           \end{array} $ \\

Further, we have to clarify in the particular case of quasi ordered classes the notion of {\it
submorphism} defined in section 2 for arbitrary categories. Given two morphisms $( A, B )$ and
$( C, D )$ in ${\cal C}$, it is easy to see that $( A, B )$ is a submorphism of $( C, D )$, if
and only if \\

(3.8) $~~~ C \leq A \leq B \leq D $ \\

Indeed, since $( A, B )$ and $( C, D )$ are morphisms in ${\cal C}$, we have $A \leq B,~ C
\leq D$. Now, $( A, B )$ being a submorphism of $( C, D )$, it follows in particular that
$( C, A )$ and $( B, D )$ are also morphisms in ${\cal C}$. Thus $C \leq A,~ B \leq D$. \\

We are now in the position to reformulate in the particular instance of quasi ordered classes
condition the (2.1) which gives the first concept of being Archimedean in the case of
arbitrary categories. Namely, in view of (3.7), (3.8), we obtain \\

(3.9) $~~~ \begin{array}{l}
       \exists~~ U \leq W ~~\mbox{objects in}~~ {\cal C} ~: \\ \\
       \forall~~ A \leq B~~\mbox{objects in}~~ {\cal C} ~: \\ \\
       \exists~~ U_1 \leq U_2 \leq \ldots \leq U_n \leq U_{n + 1 } ~~\mbox{objects in}~~
                                                  {\cal C} ~: \\ \\
       ~~~~ 1)~~ U \approx U_m \leq U_{m + 1} \approx W,~~~ 1 \leq m \leq n \\ \\
       ~~~~ 2)~~ U_1 \leq A \leq B \leq U_{n + 1}
            \end{array} $ \\

Here we note that there are two versions of this condition. Namely, if $n = 1$ in (3.9), then
we obtain \\

$~~~ \begin{array}{l}
       \exists~~ U \leq W ~~\mbox{objects in}~~ {\cal C} ~: \\ \\
       \forall~~ A \leq B~~\mbox{objects in}~~ {\cal C} ~: \\ \\
       \exists~~ U_1 \leq U_2 ~~\mbox{objects in}~~ {\cal C} ~: \\ \\
       ~~~~ 1)~~ U \approx U_1 \leq U_2 \approx W \\ \\
       ~~~~ 2)~~ U_1 \leq A \leq B \leq U_2
      \end{array} $ \\ \\

which obviously simplifies to \\

(3.10) $~~~ \begin{array}{l}
                 \exists~~ U \leq W ~~\mbox{objects in}~~ {\cal C} ~: \\ \\
                 \forall~~ A \leq B~~\mbox{objects in}~~ {\cal C} ~: \\ \\
                 ~~~~ U \leq A \leq B \leq W
             \end{array} $ \\ \\

When however $n \geq 2$ in (3.9), then quite surprisingly, we obtain the stronger condition \\

$~~~ \begin{array}{l}
       \exists~~ U \leq W ~~\mbox{objects in}~~ {\cal C} ~: \\ \\
       \forall~~ A \leq B~~\mbox{objects in}~~ {\cal C} ~: \\ \\
       \exists~~ U_1 \leq U_2 \leq \ldots \leq U_n \leq U_{n + 1 } ~~\mbox{objects in}~~
                                                  {\cal C} ~: \\ \\
       ~~~~ 1)~~ U \approx U_1 \approx \ldots \approx U_{m + 1} \approx W \\ \\
       ~~~~ 2)~~ U \approx A \approx B \approx W
            \end{array} $ \\ \\

which, when simplified, becomes \\

$~~~ \begin{array}{l}
        \exists~~ U \leq W ~~\mbox{objects in}~~ {\cal C} ~: \\ \\
        \forall~~ A \leq B~~\mbox{objects in}~~ {\cal C} ~: \\ \\
        ~~~~ U \approx A \approx B \approx W
     \end{array} $ \\ \\

thus a particular case of (3.10). In this way, we obtained \\

{\bf Proposition 1.} \\

A quasi ordered category ${\cal C}$ is Archimedean in the sense of (2.1), if and only if it
is bounded. \\

{\bf Proof.} \\
In view of (3.8), condition (3.10) obviously means that ${\cal C}$ is bounded in the sense of
the definition following (2.5). \\

{\bf Remark 1.} \\

In view of Proposition 1 above, in the case of quasi ordered categories, the concept of being
Archimedean given in (2.1) does {\it not} recover any of the two concepts (1.2) or (1.4) which
are usual in the particular case of partially ordered semigroups. \\
However, quasi ordered categories, let alone, arbitrary categories for which the concept of
being Archimedean in (2.1) was defined, have an obviously weaker structure than partially
ordered semigroups. Thus one cannot expect definition (2.1) to be able to fully include usual
concepts of being Archimedean in such a richer structure like partially ordered semigroups. \\

Let us now turn to the second concept of being Archimedean as given in (2.5), and consider it
in the particular case of quasi ordered categories ${\cal C}$. \\

In view of (3.8), a class ${\cal N}$ of morphisms in ${\cal C}$ is bounded, if and only if
there is a morphism $A \leq B$ in ${\cal C}$, such that for every morphism $C \leq D$ in
${\cal N}$, we have $A \leq C \leq D \leq B$. \\

Now in order to elucidate condition (2.5) in the case of quasi ordered categories ${\cal C}$,
we have to clarify the corresponding particular instances of (2.3). Let therefore $U \leq W$
be a morphism in ${\cal C}$. Then $\mathbb{N} ( U, W )$ is, according to (2.4), the class of
morphisms in ${\cal C}$, of the form \\

(3.11) $~~~ U_1 \leq U_{n + 1} $ \\

where \\

(3.12) $~~~ U_1 \leq \ldots \leq U_n \leq U_{n + 1} $ \\

are objects in ${\cal C}$, and in view of (3.7), they satisfy \\

(3.13) $~~~  U \approx U_m \leq U_{m + 1} \approx W,~~~ 1 \leq m \leq n $ \\

Clearly, it follows that $\mathbb{N} ( U, W )$ is always bounded, no matter which would be the
morphism $U \leq W$ in ${\cal C}$. \\

Consequently, condition (2.5) becomes \\

(3.14) $~~~ \forall~~ U \leq W ~~\mbox{morphism in}~~ {\cal C} ~:~~~ U = W $ \\

Therefore, we obtain \\

{\bf Proposition 2.} \\

A quasi ordered category ${\cal C}$ is Archimedean in the sense of (2.5), if and only if all
its morphisms are identities, that is, it is discrete, [H \& S, p. 17]. \\

{\bf Remark 2.} \\

The comments at Remark 1 above apply again. \\

\end{document}